\documentclass{amsart}
\usepackage{amsmath,amsthm,amssymb,amsfonts,amscd}
\usepackage{array}
\usepackage{comment}

\theoremstyle{plain}
\newtheorem{Thm}{Theorem}[section]
\newtheorem{Prop}[Thm]{Proposition}
\newtheorem{Coro}[Thm]{Corollary}
\newtheorem{Lem}[Thm]{Lemma}
\newtheorem{Conj}[Thm]{Conjecture}

\theoremstyle{definition}

\newtheorem{Ex}[Thm]{Example}

\newcommand{\G}{\mathcal{G}}
\newcommand{\M}{\mathcal{M}}
\newcommand{\gin}{\mathrm{gin}}
\newcommand{\Reg}{\mathrm{reg}}
\newcommand{\Codim}{\mathrm{codim}}
\newcommand{\iin}{\mathrm{in}}
\newcommand{\supp}{\mathrm{supp}}

\newcommand{\SSP}{strong Stanley property}

\newcommand{\RA}{\rightarrow}

\numberwithin{equation}{section}

\subjclass[2000]{13A02, 13C05, 13D40, 13E10, 13P10}

\begin{document}
\allowdisplaybreaks
\title[Conditions for Gin to be almost reverse lexicographic]
      {Conditions for generic initial ideals to be almost reverse lexicographic}

\author[Young Hyun Cho and Jung Pil Park]
{Young Hyun Cho$^{\dag}$ and Jung Pil Park$^{\ddag}$}
\address{$\dag$ Department of Mathematical Sciences and Research Institute for Mathematics,
                Seoul National University,
                Seoul 151-747, South Korea }
\email{youngcho@math.snu.ac.kr} %

\address{$\ddag$ National Institute for Mathematical Sciences(NIMS),
                 Daejeon 305-340, South Korea }
\email{jppark@nims.re.kr} %

\date{\today}

\thanks{The first author is partially supported by BK21, and the second author
is supported by National Institute for Mathematical Sciences(NIMS)}

\begin{abstract}
   Let $I$ be a homogeneous Artinian ideal in a polynomial
ring $R=k[x_1,\ldots,x_n]$ over a field $k$ of characteristic $0$.
We study an equivalent condition for the generic initial ideal
$\gin(I)$ with respect to reverse lexicographic order to be almost
reverse lexicographic. As a result, we show that Moreno-Socias
conjecture implies Fr\"{o}berg conjecture. And for the case $\Codim
I \le 3$, we show that $R/I$ has the strong Lefschetz property if
and only if $\gin(I)$ is almost reverse lexicographic. Finally for a
monomial complete intersection Artinian ideal
$I=(x_1^{d_1},\ldots,x_n^{d_n})$, we prove that $\gin(I)$ is almost
reverse lexicographic if $d_i > \sum_{j=1}^{i-1} d_j - i + 1$ for
each $i \ge 4$. Using this, we give a positive partial answer to
Moreno-Socias conjecture, and to Fr\"{o}berg conjecture.
\end{abstract}

\maketitle

\section{\sc Introduction}
  Let $R=k[x_1,\ldots,x_n]$ be the polynomial ring over a field
$k$. Throughout this paper, we assume that $k$ is a field of
characteristic $0$ and we use only the reverse lexicographic order
as a multiplicative term order. A monomial ideal $I$ in $R$ is said
to be {\itshape almost reverse lexicographic} if $I$ contains every
monomial $M$ which is bigger than a minimal generator of $I$ having
the same degree with $M$. One of conjectures associated with the
almost reverse lexicographic ideal is Moreno-Socias conjecture
\cite{Mo}.

\begin{Conj}[Moreno-Socias]\label{Conj:Moreno-Socias}
  If $I$ is a homogeneous ideal generated by generic forms in $R$,
then the generic initial ideal $\gin(I)$ of $I$ is almost reverse
lexicographic.
\end{Conj}

  In this paper, we study the relation between the strong Lefschetz
(or Stanley) property of standard graded Artinian $k$-algebras $R/I$
and the condition for $\gin(I)$ to be almost reverse lexicographic
(Theorem \ref{Thm:MainThm}). Let $A=R/I=\bigoplus_{i=0}^{t}A_i$ be a
standard graded Artinian $k$-algebra with $t = \max \{ i \ | \, A_i
\neq 0 \}$. Then we say that $A=R/I$ has the strong Lefschetz
property if there exists a linear form $L$ such that the
multiplication $\times L^{i} : (R/I)_d \rightarrow (R/I)_{d+i}$ is
either injective or surjective for each $d$ and $i$. And $A=R/I$ is
said to have the \SSP \ if there exists a linear form $L \in R$ such
that the multiplication map $\times L^{t-2i}:A_i \rightarrow
A_{t-i}$ is bijective for each
$i=0,1,\ldots,[t/2]$. 
Note that if such a linear form $L$ exists, then it is a generic
element in $R_1$, that is, the set of such linear forms in $R_1$
form a Zariski open set of $\mathbb{P}^{n-1}$. Note that $A$ has the
strong Stanley property if and only if it has the strong Lefschetz
property and its Hilbert function is symmetric. We will abuse
notation and refer to the strong Lefschetz or Stanley properties for
$I$ rather than for $A=R/I$.

  Stanley \cite{St} and Watanabe \cite{Wa} independently showed that
any monomial complete intersection Artinian ideal has the strong
Stanley property. And there are some results on the strong Lefschetz
property \cite{HMNW,HP,MM}. But the question whether a homogeneous
Artinian ideal $I$ generated by generic forms has the strong
Lefschetz (or Stanley) property is still open even in the case $I$
is a complete intersection.

  Another longstanding conjecture on generic algebras is
Fr\"{o}berg conjecture \cite{Fr}.
\begin{Conj}[Fr\"{o}berg]
  If $I$ is a homogeneous ideal generated by generic forms
$F_1,\ldots,F_r$ in $R$ of degrees $\deg F_i = d_i$, then the
Hilbert series $S_{R/I}(z)$ of $R/I$ is given by
\[
   S_{R/I}(z) = \left|
   \frac{\prod_{i=1}^{r}(1-z^{d_i})}{(1-z)^n}\right|.
\]
\end{Conj}

   Let $I$ be a homogeneous Artinian ideal in $R$ which has
the strong Lefschetz property. For a degree $d$ form $F \in R$ we
have the following exact sequence
\[
   0 \RA ((I:F)/I)(-d) \RA R/I(-d) \xrightarrow[\quad]{\times F} (R/I)
       \rightarrow (R/I+(F))
       \rightarrow 0.
\]
If $F$ is generic, then the Hilbert function of $R/I+(F)$ is given
by
\begin{equation*}
   H(R/I+(F),t) = \max\{H(R/I,t)-H(R/I,t-d), 0 \}.
\end{equation*}
Hence the Hilbert series of $R/I+(F)$ is given by
\[
   S_{R/I+(F)}(z) = |(1-z^{d})S_{R/I}(z)|.
\]
  Hence the study for the strong Lefschetz (or Stanley)
property of standard graded $k$-algebras defined by generic forms is
closely related with the study for Fr\"{o}berg conjecture. And, as
shown by Pardue \cite{P}, Moreno-Socias conjecture implies
Fr\"{o}berg conjecture: i.e., if Moreno-Socias conjecture is true
for any number $r$ of generic forms, then Fr\"{o}berg conjecture is
also true for any $r$, we will give another proof of this at
Corollary \ref{Cor:MorenoImplesFroeberg}.

  Recently, there are some achievements for Moreno-Socias conjecture.
For codimension 2 case, Moreno-Socias \cite{Mo}, Aguirre et
al.\cite{AOT} proved Moreno-Socias conjecture is true. And Cimpoeas
\cite{C} showed that every complete intersection Artinian ideal $I$
satisfying the strong Stanley property has the almost reverse
lexicographic $\gin(I)$ for codimension 3 case.

  In this paper, we give an equivalent condition for $\gin(I)$ to be
almost reverse lexicographic in the view point of the minimal system
of generators of $\gin(I)$ (Lemma \ref{Lem:EquivCondARL} and Theorem
\ref{Thm:MainThm}). And we generalize the result of Cimpoeas: For
any homogeneous Artinian ideal $I$ of $S=k[x_1,x_2,x_3]$, the ideal
$I$ has the strong Lefschetz property if and only if $\gin(I)$ is
almost reverse lexicographic (Proposition \ref{Prop:Codim3Case}).
Then we show that $\gin(I)$ is almost reverse lexicographic for a
monomial complete intersection ideal
$I=(x_1^{d_1},\ldots,x_n^{d_n})$ in $R=k[x_1,\ldots,x_n]$ if $d_i >
\sum_{j=1}^{i-1} d_j - i + 1$ for $i \ge 4$ (Lemma
\ref{Lem:MonomialCase}). At last we show that Moreno-Socias
conjecture is true if $I$ is a complete intersection Artinian ideal
of $R$ which is generated by generic forms of degrees $d_i$
satisfying the same condition (Corollary
\ref{Cor:Moreno-SociasConjecture}). And we give a partial solution
on Fr\"{o}berg conjecture (Corollary \ref{Cor:FroebergConjecture}).
As recent works Harima et al.\cite{HW} showed similar results with
ours independently.

  Suppose that $I$ is a homogeneous Artinian ideal
of $R=k[x_1,\ldots,x_n]$. In the paper \cite{CCP}, Cho and the
authors of this paper showed that the minimal system of generators
of $\gin(I)$ is completely determined by the positive integer $f_1$
and functions $f_i:\mathbb{Z}^{i-1}_{\ge 0} \longrightarrow
\mathbb{Z}_{\ge 0} \cup \{\infty\}$ defined as follows:
   \begin{align}\label{f's}
   \begin{split}
   f_1=&\min \{t \ | \, x^t_1 \in \gin(I) \} \text{, and }\\
   f_i(\alpha_1, \ldots, \alpha_{i-1}) =& \min \{ t \ | \,
   x^{\alpha_1}_1 \cdots x^{\alpha_{i-1}}_{i-1} x^{t}_i \in
   \gin(I) \},
   \end{split}
  \end{align}
 for each $2 \le i \le n$.

\begin{Prop}\cite{CCP}\label{Prop:MinimalGen}
   Let $I$ be a homogeneous Artinian ideal in
$R=k[x_1,\ldots,x_n]$. Suppose that $f_1,\ldots,f_{n}$ are defined
for $\gin(I)$ as in $(\ref{f's})$. Then the minimal system of
generators $\G(\gin(I))$ of $\gin(I)$ is
\begin{equation*}
\begin{split}
   \G(\gin(I)) =\{x_1^{f_1}\} \cup
       \left\{
              x^{\alpha_1}_1 \cdots x^{\alpha_{i-1}}_{i-1}
              x^{f_i(\alpha_1, \ldots, \alpha_{i-1})}_i
              \left| \begin{array}{l}
                       2 \le i \le n, \\
                       0 \le \alpha_1 < f_1, \ \text{and} \\
                       0 \le \alpha_{j} < f_{j}(\alpha_1, \ldots,\alpha_{j-1}) \\
                     \text{for each} \ 2 \le j \le i
                     \end{array}
              \right.
            \right\}.
\end{split}
\end{equation*} \qed
\end{Prop}

For each $ 1 \le i \le n-1$, let the set $J_i$ be defined as
\begin{equation} \label{Set:IndexMinGen}
\begin{split}
  J_i = \left\{(\alpha_1,\ldots,\alpha_{i})
              \left| \begin{array}{l}
                       0 \le \alpha_1 < f_1, \ \text{and} \\
                       0 \le \alpha_{j} < f_{j}(\alpha_1, \ldots, \alpha_{j-1}) \\
                       \text{for each} \ 2 \le j \le i
                     \end{array}
              \right.
        \right\}.
\end{split}
\end{equation}

  We mainly use the sets $J_{i}$ to
prove the main theorem and its corollaries. Hence we need to
investigate the sets $J_{i}$ closely. In what follows, we use the
following notations for simplicity. For
$\alpha=(\alpha_1,\ldots,\alpha_i) \in \mathbb{Z}^{i}_{\ge 0}$, we
denote $\sum_{j=0}^{i} \alpha_j$ by $|\alpha|$. And for
$\beta=(\beta_1,\ldots,\beta_i) \in \mathbb{Z}^{i}_{\ge 0}$, we say
that $\beta > \alpha$ if $x_1^{\beta_1} \cdots x_i^{\beta_i} >
x_1^{\alpha_1} \cdots x_i^{\alpha_i}$.

\begin{Lem}\label{Lem:PropertyOfIndexSet}
  Let $I$ be a homogeneous Artinian ideal in $R=k[x_1,\ldots,x_n]$.
For each $1 \le i \le n-1$, let $J_{i}$ be the set defined for
$\gin(I)$ as in \rm{(}\ref{Set:IndexMinGen}\rm{)}.
\begin{enumerate}
  \item An element $(\alpha_1,\ldots,\alpha_{i}) \in \mathbb{Z}^{i}_{\ge 0}$
        belongs to $J_{i}$ if and only if
        the monomial $x_1^{\alpha_1}\cdots x_{i}^{\alpha_{i}}$ is not
        contained in $\gin(I)$.
  \item If $\alpha=(\alpha_1,\ldots,\alpha_{i}) \in J_{i}$, then
        the element $(0,\ldots,0,|\alpha|) \in \mathbb{Z}^{i}_{\ge 0}$
        belongs to $J_{i}$.
        Furthermore, we have
        $f_{i+1}(\alpha_1,\ldots,\alpha_{i}) \le f_{i+1}(0,\ldots,0,|\alpha|)$.
  \item For two elements $(0,\ldots,0,a), (0,\ldots,0,b)$ of $J_{i}$,
        if $a \le b$, then we have
        \[
           a + f_{i+1}(0,\ldots,0,a) \ge b + f_{i+1}(0,\ldots,0,b).
        \]
\end{enumerate}
\end{Lem}
\begin{proof}
  \begin{enumerate}
  \item The assertion follows easily from Proposition \ref{Prop:MinimalGen}
        and the definition of $J_i$.
  \item If $(0,\ldots,0,|\alpha|)$ is not an element of $J_{i}$,
        then the monomial $x_{i}^{|\alpha|}$ belongs to $\gin(I)$ as shown in (1).
        But this implies that the monomial $x_1^{\alpha_1}\cdots x_{i}^{\alpha_{i}}$ is
        also contained in $\gin(I)$, since $\gin(I)$ is strongly stable.
        This contradicts to $(\alpha_1,\ldots,\alpha_{i}) \in J_{i}$.
        And the last assertion follows from strongly stableness of $\gin(I)$ and the definition of $f_{i+1}$.
  \item Set $\mu = f_{i+1}(0,\ldots,0,a)$ and $t = \min\{b-a,\mu\}$.
        By the definition of $f_{i+1}$, we have $x_{i}^a x_{i+1}^{\mu} \in \gin(I)$.
        And since $\gin(I)$ is strongly stable, this implies that
        $x_{i}^{a+t} x_{i+1}^{\mu-t} \in \gin(I)$. If $\mu \le b-a$,
        then $x_i^{a+t}x_{i+1}^{\mu-t}=x_i^{a+\mu} \in \gin(I)$, and hence
        $x_i^{b}$ is also contained in $\gin(I)$. But, this contradicts
        to $(0,\ldots,0,b) \in J_i$. Hence $\mu > b-a$ and
        $x_i^b x_{i+1}^{\mu-t} = x_{i}^{a+t} x_{i+1}^{\mu-t} \in
        \gin(I)$. This shows that
        \[
           f_{i+1}(0,\ldots,0,b) \le \mu-t = f_{i+1}(0,\ldots,0,a) - (b-a).
        \]
  \end{enumerate}
\end{proof}

   In the paper \cite{ACP}, Ahn et al. gave the following tool
detecting whether $I$ has the strong Lefschetz (or Stanley)
property, from the view point of the minimal system of generators of
$\gin(I)$. This tool gives us the chance to prove the main theorems
easily.

\begin{Prop}\label{Prop:SLPAndSSP}\cite{ACP}
   Let $I$ be a homogeneous Artinian ideal in $R=k[x_1,\ldots,x_n]$
with $t = \max\{i | \, (R/I)_i \neq 0 \}.$
\begin{enumerate}
  \item $I$ has the strong Lefschetz property if and only if we have
        \[
            f_n(0,\ldots,0,|\alpha| + 1) + 1  \le f_n(\alpha),
        \]
        for any $\alpha=(\alpha_1,\ldots,\alpha_{n-1}) \in J_{n-1}$.
  \item $I$ has the strong Stanley property if and only if we have
        \[
            f_n(\alpha) = t - 2|\alpha| + 1,
        \]
        for any $\alpha=(\alpha_1,\ldots,\alpha_{n-1}) \in J_{n-1}$.
\end{enumerate}
\end{Prop}

\vskip 1cm

\section{\sc Main Results}
  Let $I$ be a homogeneous Artinian ideal in $R=k[x_1,\ldots,x_n]$.
Unless otherwise stated, we assume that the integer $f_1$, the
functions $f_i$ and the sets $J_i$ are defined for $\gin(I)$ as in
(\ref{f's}) and (\ref{Set:IndexMinGen}), respectively.

The following lemma gives an equivalent condition for $\gin(I)$ to
be almost reverse lexicographic.

\begin{Lem}\label{Lem:EquivCondARL}
  Let $I$ be a homogeneous Artinian ideal in
$R=k[x_1,\ldots,x_n]$. Then $\gin(I)$ is almost reverse
lexicographic if and only if for any $1 \le i \le n-1$ the following
conditions are satisfied:
\begin{enumerate}
\item $f_{i+1}(0,\ldots,0,|\alpha|+1)+1 \le f_{i+1}(\alpha)$
      for any $\alpha = (\alpha_1,\ldots,\alpha_i) \in J_i$, and
\item $f_{i+1}(\beta) \le f_{i+1}(\alpha)$ for any $\alpha = (\alpha_1,\ldots,\alpha_i), \beta
=(\beta_1,\ldots,\beta_i) \in J_i$ with $|\alpha| = |\beta|$ and
$\alpha < \beta$.
\end{enumerate}
\end{Lem}
\begin{proof}
$(\Rightarrow):$  Fix $i$. Note that $f_{i+1}(\alpha) > 0$,
otherwise $x_1^{\alpha_1} \cdots x_i^{\alpha_i} \in \gin(I)$ by the
definition of $f_{i+1}$, and this contradicts to $\alpha \in J_{i}$
as shown in Lemma \ref{Lem:PropertyOfIndexSet} (1). Hence we have
\[
   x_{i}^{|\alpha|+1} x_{i+1}^{f_{i+1}(\alpha)-1}
   > x_1^{\alpha_1}\cdots x_{i}^{\alpha_i} x_{i+1}^{f_{i+1}(\alpha)}.
\]
This shows that the monomial $x_{i}^{|\alpha|+1}
x_{i+1}^{f_{i+1}(\alpha)-1}$ is contained in $\gin(I)$, since
$\gin(I)$ is almost reverse lexicographic.
So the assertion (1) follows by the definition of $f_{i+1}$. And the
assertion (2) follows easily by the hypothesis and the definition of
$f_{i+1}$.

$(\Leftarrow):$ Let $1 \le i \le n-1$, and let $M = x_1^{\beta_1}
\cdots x_{i-1}^{\beta_{i-1}}x_i^{b}$, $N = x_1^{\alpha_1} \cdots
x_{i-1}^{\alpha_{i-1}}x_i^{a}$ be  monomials in $R$ having the same
degree. Suppose that $M > N$ and $N$ is a minimal generator of
$\gin(I)$. If we set $\alpha=(\alpha_1,\ldots,\alpha_{i-1}),
\beta=(\beta_1,\ldots,\beta_{i-1})$, then we have $\alpha \in
J_{i-1}$, $a = f_i(\alpha)$ and $b \le a$. We have to show that $M$
belongs to $\gin(I)$. We may assume $\beta \in J_{i-1}$ by Lemma
\ref{Lem:PropertyOfIndexSet} (1). Hence we have to show that
$f_i(\beta) \le b$.

If $b = a$, then $|\alpha| = |\beta|$ and $x_1^{\beta_1}\cdots
x_{i-1}^{\beta_{i-1}} > x_1^{\alpha_1} \cdots
x_{i-1}^{\alpha_{i-1}}$. By the condition (2), we have
\[
  f_{i}(\beta) \le f_{i}(\alpha) = a = b.
\]

If $b < a$, then $|\beta| > |\alpha|$. Hence we have
\begin{align*}
   |\beta| + f_{i}(\beta) & \le |\beta| + f_{i}(0,\ldots,0, |\beta|) \\
                          & \le |\alpha| + 1 + f_{i}(0,\ldots,0,|\alpha|+1) \\
                          & \le |\alpha| + f_{i}(\alpha)
                          = |\alpha| + a = |\beta| + b,
\end{align*}
where the first and second inequalities follow from Lemma
\ref{Lem:PropertyOfIndexSet} (2) and (3), respectively, and the
third one follows from the condition (1). This shows that
$f_{i}(\beta) \le b$.
\end{proof}

\begin{Coro}\label{Cor:Codim2Case}
 For every homogeneous Artinian ideal $K$ in the polynomial ring
$S=k[x_1,x_2]$, $\gin(K)$ is almost reverse lexicographic.
\end{Coro}
\begin{proof}
   If $\alpha \in J_1$, then $x_1^{\alpha}x_2^{f_2(\alpha)} \in \gin(K)$ and $f_2(\alpha) > 0$.
   Since $\gin(K)$ is strongly stable, we have
   $x_1^{\alpha+1}x_2^{f_2(\alpha)-1} \in \gin(K)$. Hence
   $f_2(\alpha+1) \le f_2(\alpha)-1$. This shows that two
   conditions (1) and (2) in Lemma \ref{Lem:EquivCondARL} are
   satisfied.
\end{proof}

    Note that we must verify that two conditions (1) and (2)
in Lemma \ref{Lem:EquivCondARL} are satisfied for every $i$ from 1
to $n-1$, in order to show that $\gin(I)$ is almost reverse
lexicographic. Using Corollary \ref{Cor:Codim2Case}, we can reduce
the range of $i$ to check, furthermore we can rephrase Lemma
\ref{Lem:EquivCondARL} as follows.

\begin{Prop}\label{Prop:SecondStep}
  Let $I$ be a homogeneous Artinian ideal in the polynomial ring
$R=k[x_1,\ldots,x_n]$. Then $\gin(I)$ is almost reverse
lexicographic if and only if for each $2 \le i \le n-1$ the
following conditions are satisfied:
\begin{enumerate}
   \item There exist generic linear forms
         $L_1,\ldots,L_{n-i-1}$ in $R$ such that the ring
         $R/I+(L_1,\ldots,L_{n-i-1})$ has the strong Lefschetz property.
   \item For any  two elements $\alpha = (\alpha_1,\ldots,\alpha_i)$
         and $\beta =(\beta_1,\ldots,\beta_i)$  of $J_i$
         with $|\alpha| = |\beta|$ and $\alpha < \beta$, we have
         $f_{i+1}(\beta) \le f_{i+1}(\alpha)$.
\end{enumerate}
\end{Prop}
\begin{proof}
   Let $1 \le i \le n-1$. For generic linear forms $L_1,\ldots,L_{n-i-1}$
in $R$, if we set $\overline{I}$ to be the image of $I$ in the ring
$R/(L_1,\ldots,L_{n-i-1})$, then $\gin(\overline{I}) =
\gin(I)_{x_n\rightarrow 0, \ldots, x_{i+2} \rightarrow 0}$ by Green
(see Corollary 2.15 in the paper \cite{Gr}, and see also Wiebe
\cite{Wie}). If we denote $f^I_j$ and $J^I_j$ the functions and the
sets defined for $\gin(I)$ as in (\ref{f's}) and
(\ref{Set:IndexMinGen}) respectively, then this implies that
$f^I_{j+1} = f^{\overline{I}}_{j+1}$ and $J^I_j= J^{\overline{I}}_j$
for any $1 \le j \le i$. Since
\[
   \frac{R/(L_1,\ldots,L_{n-i-1})}{\overline{I}} =
   R/I+(L_1,\ldots,L_{n-i-1}),
\]
the first condition in Lemma \ref{Lem:EquivCondARL} implies that
$R/I+(L_1,\ldots,L_{n-i-1})$ has the strong Lefschetz property for
each $1 \le i \le n-1$. But since every homogeneous Artinian ideal
$K$ of codimension $2$ has almost reverse lexicographic $\gin(K)$,
it is enough to check only for $i$ from $2$ to $n-1$.

\end{proof}

 As shown in Introduction, Proposition \ref{Prop:SecondStep} implies that

\begin{Coro}\label{Cor:MorenoImplesFroeberg}
   Moreno-Socias conjecture implies that Fr\"{o}berg conjecture, that is,
if Moreno-Socias conjecture is true for any number $r$ of generic
forms in a polynomial ring $R=k[x_1,\ldots,x_n]$, then Fr\"{o}berg
conjecture is also true for any $r$.
\end{Coro}

   The following example shows the second condition in Proposition \ref{Prop:SecondStep}
cannot be omitted.
\begin{Ex} \label{Ex:DoubleSLPnotImplyALR}
  Consider the following strongly stable monomial ideal
  \[ I = \left(
         \begin{array}{l}
         x^2, xy^2, y^4, y^3z, xyz^2, y^2z^2, xz^3, yz^3, z^4, \\
         y^3w, xyzw, \underline{xz^2w^2, y^2zw^3}, yz^2w^3, z^3w^3, \\
         xyw^4, y^2w^4, xzw^4, yzw^4, z^2w^4 \\
         xw^5, yw^5, zw^5, w^6
         \end{array}
         \right) \subset S=k[x,y,z,w].
  \]
  Note that both $S/I$ and $S/I+(L)$ have the strong Lefschetz property
by Proposition \ref{Prop:SLPAndSSP}. Although $xz^2w^2$ is a minimal
generator of $I$ and $y^2zw^2 > xz^2w^2$, $y^2zw^2$ does not belong
to $I$. Hence $I$ is not almost reverse lexicographic.
\end{Ex}

  But if we restrict our interests to the case that $I$ is
a homogeneous Artinian ideal in $S=k[x_1,x_2,x_3]$, then we can show
that the second condition in Proposition \ref{Prop:SecondStep} is
superfluous, that is, $\gin(I)$ is almost reverse lexicographic if
and only if $S/I$ has the strong Lefschetz property.

\begin{Prop}\label{Prop:Codim3Case}
  Let $I$ be a homogeneous Artinian ideal of $S=k[x_1,x_2,x_3]$.
Then $S/I$ has the strong Lefschetz property if and only if
$\gin(I)$ is almost reverse lexicographic.
\end{Prop}
\begin{proof}
$(\Rightarrow)$  It suffices to show that the second condition in
Proposition \ref{Prop:SecondStep} is fulfilled. Let
$\alpha=(\alpha_1,\alpha_2)$, $\beta=(\beta_1,\beta_2) \in J_2$. If
$|\alpha| = |\beta|$ and $\alpha < \beta$, then $\alpha_2 \ge
\beta_2$ and $\beta_1 = \alpha_1 + (\alpha_2 - \beta_2)$. Since the
monomial $x_1^{\alpha_1}x_2^{\alpha_2}x_3^{f_3(\alpha)}$ is
contained in $\gin(I)$, and since $\gin(I)$ is strongly stable, we
have
\[
  x_1^{\beta_1}x_2^{\beta_2}x_3^{f_3(\alpha)} =
  x_1^{\alpha_1 + (\alpha_2 - \beta_2)} x_2^{\alpha_2 - (\alpha_2 - \beta_2)}
  x_3^{f_3(\alpha)} \in
  \gin(I).
\]
This shows that $f_3(\beta) \le f_3(\alpha)$ by the definition of
$f_3$.

$(\Leftarrow)$ This follows from Proposition \ref{Prop:SecondStep}.
\end{proof}

  Putting Propositions \ref{Prop:SecondStep} and \ref{Prop:Codim3Case} together,
we obtain the main result of this paper.

\begin{Thm}\label{Thm:MainThm}
  Let $I$ be a homogeneous Artinian ideal in the polynomial ring
$R=k[x_1,\ldots,x_n]$. Then $\gin(I)$ is almost reverse
lexicographic if and only if the following conditions are satisfied:
\begin{enumerate}
   \item For any $0 \le i \le n-3$, there exist generic linear forms
         $L_1,\ldots,L_{i}$ in $R$ such that the ring
         $R/I+(L_1,\ldots,L_{i})$ has the strong Lefschetz property.
   \item For each $3 \le i \le n-1$, if $\alpha = (\alpha_1,\ldots,\alpha_i)$,
         $\beta =(\beta_1,\ldots,\beta_i)$ are two elements of $J_i$
         with $|\alpha| = |\beta|$ and $\alpha < \beta$, then we have
         $f_{i+1}(\beta) \le f_{i+1}(\alpha)$.
\end{enumerate}
\end{Thm}

  At Example \ref{Ex:DoubleSLPnotImplyALR}, we showed that the first condition
in Proposition \ref{Prop:SecondStep} is not enough to ensure that
$\gin(I)$ is almost reverse lexicographic, if the number of
variables of the ring $R$ is greater than or equal to 4. But for any
$0 \le i \le n-3$, if there exist generic linear forms $L_1, \ldots,
L_i$ such that $I+(L_1,\ldots,L_i)$ has the strong Stanley property,
then $\gin(I)$ is almost reverse lexicographic as shown in the
following corollary.

\begin{Coro}\label{Cor:MoreVariableWithSSP}
Let $I$ be a homogeneous Artinian ideal in $R=k[x_1,\ldots,x_n]$.
For each $0 \le i \le n-3$, if there exist generic linear forms
$L_1,\ldots,L_i$ in $R$ such that $I+(L_1,\ldots,L_i)$ has the
strong Stanley property, then $\gin(I)$ is almost reverse
lexicographic.
\end{Coro}
\begin{proof}
  It is enough to show that the second condition in Theorem \ref{Thm:MainThm} is
satisfied for each $3 \le i \le n-1$. Fix $i$ and let
$\alpha=(\alpha_1,\ldots,\alpha_i)$ and
$\beta=(\beta_1,\ldots,\beta_i)$ be two elements of $J_i$ with
$|\alpha| = |\beta|$ and $\beta
> \alpha$. By the assumption, there exist generic linear forms
$L_1,\ldots,L_{n-i-1}$ in $R$ such that $I+(L_1,\ldots,L_{n-i-1})$
has the strong Stanley property. From the reason described in the
proof of Proposition \ref{Prop:SecondStep}, we have
\[
   f_{i+1}(\beta) = t - 2|\beta| + 1 = t - 2|\alpha| + 1 =
   f_{i+1}(\alpha),
\]
where $t = \max\{j |\, (R/I+(L_1,\ldots,L_{n-i-1}))_j \neq 0 \}$.
Hence the assertion follows.
\end{proof}

  As a result, we will show that Moreno-Socias conjecture is true
for the case that $K$ is a complete intersection Artinian ideal
generated by generic forms $F_1,\ldots,F_n$ in $R$ of degrees $d_i$
with $d_i > \sum_{j=1}^{i-1} d_j - i +1$ for each $i \ge 4$. To do
so, we will show first that if $I=(x_1^{d_1},\ldots,x_n^{d_n})$
under the same condition on the $d_i$, then $\gin(I)$ is almost
reverse lexicographic. We use the result of Stanley and Watanabe:
every monomial complete intersection Artinian ideal has the strong
Stanley property (see \cite{St} or Corollary 3.5 \cite{Wa}).

\begin{Lem}\label{Lem:MonomialCase}
Let $I=(x_1^{d_1},\ldots,x_n^{d_n})$. If $d_i > \sum_{j=1}^{i-1} d_j
- i + 1$ for each $i \ge 4$, then $\gin(I)$ is almost reverse
lexicographic.
\end{Lem}
\begin{proof}
  By Corollary \ref{Cor:MoreVariableWithSSP}, it suffices to show
that for each $0 \le i \le n-3$, there exist generic linear forms
$L_1,\ldots,L_i$ such that $R/I+(L_1,\ldots,L_i)$ has the strong
Stanley property. Hence it is enough to show that
$R/I+(L_1,\ldots,L_i)$ is isomorphic to
$k[x_1,\ldots,x_{n-i}]/(x_1^{d_1},\ldots,x_{n-i}^{d_{n-i}})$ for
each $0 \le i \le n-3$. We will show this by induction on $i$. At
first, for simplicity, we denote by $\M_j$ the set of monomials
$\{x_1^{d_1},\ldots,x_{n-j}^{d_{n-j}}\}$ for each $0 \le j \le n-3$.

  For the case $i=0$, the claim is true by the hypothesis.
Assume that the claim is true for the case $i = s < n-3$. Then there
exist generic linear forms $L_1,\ldots,L_s$ such that
$R/I+(L_1,\ldots,L_s)$ is isomorphic to $S/(\M_s)$, where
$S=k[x_1,\ldots,x_{n-s}]$. Choose a generic linear form $L_{s+1} \in
R$ such that the ideal generated by $\M_{s+1} \cup
\{L_{1},\ldots,L_{s},L_{s+1}\}$ is a complete intersection Artinian
ideal in $R$. If we set $L^{'}$ to be the image of $L_{s+1}$ in
$S/(\M_s)$, then we have $R/I+(L_1,\ldots,L_{s+1}) =
S/(\M_s)+(L^{'})$. Let $K$ be the ideal generated by $\M_{s-1} \cup
\{L^{''}\}$ in $S$, where $L^{''}$ is a generic form in $S$ such
that the image of $L^{''}$ in $S/(\M_{s})$ is $L^{'}$. Note that $K$
is a complete intersection Artinian ideal in $S$, and that the
Castelnuovo-Mumford regularity of the ideal $K$ is
\[
   \Reg \, K = 1 + \sum_{j=1}^{n-s-1} d_{j} - (n-s-1).
\]
Since $s < n-3$, we have $d_{n-s} \ge \Reg \, K$. This implies that
$x_{n-s}^{d_{n-s}} = 0$ in $S/K$. Hence we have
\begin{align*}
   R/I+(L_1,\ldots,L_{s},L_{s+1})
      &= S/(\M_s)+(L^{'})
       = S/K+(x_{n-s}^{d_{n-s}}) = S/K \\
      &= k[x_1,\ldots,x_{n-s-1},x_{n-s}]/(\M_{s-1}) + (L^{''}) \\
      &= k[x_1,\ldots,x_{n-s-1}]/(\M_{s-1}),
\end{align*}
the last equation follows since $L^{''}$ is generic. So we are done.
\end{proof}

  And we need to know another way to compute generic initial ideals.
The following is introduced in Eisenbud \cite{E}:

  By a monomial of $\wedge^d R_i$ we mean an element of the form
$N=n_1 \wedge \cdots \wedge n_d$, where the $n_j$ are degree $i$
monomials of $R$, and we denote the support of $N$ by
$\supp(N)=\{n_1,\ldots,n_d\}$. We define a term in $\wedge^d R_i$ to
be a product $p \cdot N$, where $p \in k$ and $N$ is a monomial in
$\wedge^d R_i$. We will say that $p \cdot N = p \cdot n_1 \wedge
\cdots \wedge n_d$ is a normal expression if the $n_j$ are ordered
so that $n_1 > \cdots > n_d$.

  We order the monomials of $\wedge^d R_i$ by ordering their normal
expressions lexicographically: if $N=n_1 \wedge \cdots \wedge n_d$
and $N' = n'_1 \wedge \cdots \wedge n'_d$ are normal expressions,
than $N > N'$ if and only if $n_j > n'_j$ for the smallest $j$ with
$n_j \neq n'_j$. We extend the order to terms, and define the
initial term of an element $F \in \wedge^d R_i$ to be the greatest
term with respect to the given order.

  If $K$ is a homogeneous ideal of $R$, then there is a Zariski
open set $U \subset \mathrm{GL}(n,k)$ such that for each $i \ge 0$,
$\wedge^d \gin(K)_i$ is spanned by the greatest monomial of
$\wedge^d R_i$ that appears in any $\wedge^d(gK_i)$ with $g \in U$,
where $d = \dim_k K_i$, that is, if $N=n_1 \wedge \ldots \wedge n_d$
is the greatest monomial that appears in any $\wedge^d(gK_i)$ with
$g \in U$, then $\gin(K)_i = kn_1 \oplus \cdots \oplus kn_d$ (see
Theorem 15.18 in Eisenbud \cite{E} for the details).

 Hence, in order to compute $\gin(K)_i$, choose the basis
$F_1,\ldots,F_d$ of $gK_i$ with $g \in U$, where $d = \dim_k K$. Let
$N_1, \ldots, N_s$ be the monomial basis of $\wedge^d R_i$ with $N_1
> \cdots > N_s$ with respect to the given order. If we write $F_1
\wedge \cdots \wedge F_d = \sum p_j \cdot N_j$ with $p_j \in k$,
then the greatest monomial of $\wedge^d(gK_i)$ is the first monomial
$N_j$ such that $p_j \neq 0$. Note that each $p_j$ are given as
polynomial expressions in the coefficients of the $F_i$. Hence each
$p_j$ are given as polynomial expressions in the coefficients of a
minimal system of generators of $gK$. And note that if $m_1, \ldots,
m_t$ are the monomials of degree $i$ contained in $\gin(K)_{\le
i-1}$, the ideal generated by the elements of $\gin(K)$ of degrees
$\le i-1$, then $p_j = 0$ for any $j$ with $\{m_1,\ldots,m_t\}
\nsubseteq \supp(N_j)$. Hence if $M_1,\ldots,M_u$ is the monomial
basis of $\wedge^d R_i$ such that $\{m_1,\ldots,m_t\} \subset
\supp(M_j)$ and $M_1 > \cdots
> M_u$, then we can write $F_1\wedge \cdots \wedge F_d = \sum p_j
\cdot M_j$, where each $p_j \in k$ is a polynomial expression in the
coefficients of a minimal system of generators of $gK$. Note the
greatest monomial of $\wedge^d(gK_i)$ is the first monomial $M_j$
such that $p_j \neq 0$ and that the monomials in $\supp(M_j)$ except
$m_1,\ldots,m_t$ are the minimal generators of $\gin(K)$ having
degree $i$. The following theorem will give a positive answer to
Moreno-Socias conjecture in our case.

\begin{Thm}\label{Thm:SuffCondForGenIdealHavingARLGin}
  Let $I$ be a homogeneous ideal in $R=k[x_1,\ldots,x_n]$ which has
a almost reverse lexicographic generic initial ideal. If $K$ is a
homogeneous ideal generated by generic forms in $R$ such that
$H(R/I,d)=H(R/K,d)$ for all $d$, then $\gin(K)$ is also almost
reverse lexicographic.
\end{Thm}
\begin{proof}
  We will prove this by showing that $\gin(I)_i = \gin(K)_i$ for
every $i \ge 0$. It is clear for the case $i=0$.

Suppose that $\gin(I)_j = \gin(K)_j$ for $j=0,\ldots,i$. Let
$d=\dim_k I_{i+1} = \dim_k K_{i+1}$. Without loss of generality, we
may assume that we make a general choice of coordinates for $K$ and
$I$, i.e. $\gin(K) = \iin(K)$ and $\gin(I) = \iin(I)$. Let
$m_1,\ldots,m_t$ be the monomials of degree $i+1$ in $\gin(I)_{\le
i} = \gin(K)_{\le i}$. And let $n_{t+1},\ldots,n_{d}$ be the minimal
generators of $\gin(I)$ of degree $i+1$ with $n_{t+1}>\cdots>n_{d}$.
Then we know that the monomial $N=m_1\wedge \cdots \wedge m_t \wedge
n_{t+1} \wedge \cdots \wedge n_d$ is the greatest monomial appearing
in $\wedge^d I_{i+1}$. Let $M_1,\ldots,M_u$ be the monomial basis of
$\wedge^d R_{i+1}$ such that $\{m_1,\ldots,m_t\} \subset \supp(M_j)$
and $M_1 > \cdots > M_u$. If $F_1,\ldots,F_d$ is a basis of
$K_{i+1}$, then $F_1\wedge\cdots\wedge F_d$ can be written as
$F_1\wedge\cdots\wedge F_d = \sum p_j \cdot M_j$ with $p_j \in k$.

  Now we will show that $N=M_1$. Let $M_1 = m_1\wedge\cdots\wedge
m_t \wedge l_{t+1}\ldots \wedge l_{d}$ for some monomials
$l_{t+1},\ldots,l_{d}$ in $R_{i+1}$ with $l_{t+1} > \cdots > l_d$.
It suffices to show that $n_{t+1}\wedge\cdots\wedge n_{d} \ge
l_{t+1}\wedge\cdots\wedge l_d$ as monomials in $\wedge^{d-t}
R_{i+1}$ by the choice of $M_1$. But this is clear because $\gin(I)$
is almost reverse lexicographic and $n_{t+1},\ldots,n_{d}$ are
minimal generators of $\gin(I)$ having degrees $i+1$.

  This show that $p_1$, the coefficient of $M_1$, is a nonzero
polynomial expression in the coefficients of a minimal system of
generators of $K$. Since $K$ is generated by generic forms and $k$
is a field of characteristic $0$, this implies that $p_1 \neq 0$.
Hence we have $\gin(I)_{i+1} = \{m_1,\ldots,m_t, n_{t+1}, \ldots,
n_d \} = \gin(K)_{i+1}$.

  Since $k$ is a field of characteristic 0 and $\gin(K)$ is finitely
generated, for each step $i$, we can have only finitely many closed
subsets in a projective space which are defined by polynomials
having the variables as the coefficients of a minimal system of
generators of $K$. Since $K$ is generated by generic forms, we have
$\gin(I)_i = \gin(K)_i$ for all $i \ge 0$.
\end{proof}

\begin{Coro}\label{Cor:Moreno-SociasConjecture}
  Let $K$ be a generic complete intersection Artinian ideal
generated by generic forms $G_1,\ldots,G_n$ in $R=k[x_1,\ldots,x_n]$
with $\deg G_i = d_i$, if $d_i > \sum_{j=1}^{i-1} d_j - i + 1$ for
each $i \ge 4$, then $\gin(K)$ is almost reverse lexicographic. In
particular, every generic complete intersection Artinian ideal of
codimension 3 has almost reverse lexicographic generic initial
ideal, and hence has the strong Stanley property.
\end{Coro}
\begin{proof}
Let $I=(x_1^{d_1},\ldots,x_n^{d_n})$. Then $H(R/I,d)=H(R/K,d)$ for
all $d$. By Theorem \ref{Thm:SuffCondForGenIdealHavingARLGin}, we
are done.
\end{proof}

This gives a positive partial answer for Fr\"{o}berg conjecture.

\begin{Coro}\label{Cor:FroebergConjecture}
  Let $K$ be a generic complete intersection Artinian ideal generated
by generic forms $F_1, \ldots, F_n$ in $R$ of degrees $d_i$ with
$d_i > \sum_{j=1}^{i-1} d_j - i + 1$ for each $i \ge 4$. If $F$ is a
generic form of degree $d$, then the Hilbert series of $R/I+(F)$ is
\[
   S_{R/I+(F)}(z) = \left|
                    \frac{(1-z^d)\prod_{i=1}^{n}(1-z^{d_i})}{(1-z)^n}\right|.
\]
\end{Coro}

\vskip 1cm
\bibliographystyle{amsalpha}

\end{document}